\documentclass{amsart}
\usepackage{amsmath}%
\usepackage{amsfonts}%
\usepackage{amssymb}%
\usepackage{graphicx}
\theoremstyle{plain}
\newtheorem{acknowledgement}{Acknowledgement}

\newtheorem{conclusion}{Conclusion}

\newtheorem{definition}{Definition}

\numberwithin{equation}{section}
\begin{document}
\title[Probability distribution of distances between local extrema]{Probability distribution of distances between local extrema of random number series}
\author{Argyn Kuketayev}
\email[A. Kuketayev]{jawabean@gwu.edu}%
\date{August 12, 2006}

\begin{abstract}
There is a sequence of random numbers $x_1,x_2 \ldots,x_{n}$ and so on. Numbers are independent of each other, but all numbers are from the same continuous distribution. If $ x_{i-1} < x_i > x_{i+1}$, then $x_i$ is a local maximum. Here, we show that the probability mass function (PMF) $f_m(d)$ of distances $d$ between local maxima is non-parametric and the same for any probability distribution of random numbers in the sequence, and that the average distance is exactly 3. We present a method of computation of this PMF and its table for distances betwen 2 and 29. This PMF is confirmed to match distance distributions of sample random number sequences, which were created by pseudo-random number generators or obtained from "true" random number sources.
\end{abstract}
\maketitle

\section{Average distnace between local maxima}

\begin{figure}[h]
	\centering
		\includegraphics{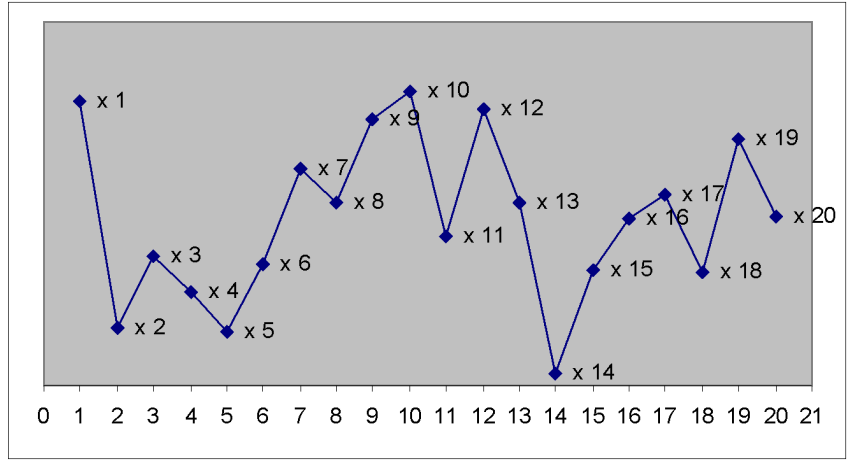}
	\caption{Sample random sequence}
	\label{sample-seq}
\end{figure}

\noindent Let's take any number in the sequence and find out the probability that it's a local maximum. 
\begin{definition}
A number $x_i$ is a local maximum, if the following condition is true $ x_{i-1} < x_i > x_{i+1}$.
\end{definition}
First, we'll use a combinatorial approach. Consider the following sequence (Figure \ref{sample-seq}) of pseudo-random numbers generated by MS Excel RAND() function:
\begin{enumerate}
	\item 0.935536495
\item 0.191531578
\item 0.429049655
\item 0.308968021
\item 0.179540986
\item 0.401329789
\item 0.71581906
\item 0.604617962
\item 0.877254876
\item 0.973280207
\item 0.489033299
\item 0.912367351
\item 0.604552972
\item 0.039395302
\item 0.3780448
\item 0.55317569
\item 0.6308772
\item 0.373163479
\item 0.812434426
\item 0.560173882
\end{enumerate}

Let's take any three consequitive numbers $x_{i-1}$, $x_i$ and $x_{i+1}$. For instance, for $i=3$ we have: $x_2=0.191531578$, $x_3=0.429049655$ and $x_4=0.308968021$. In this case, $x_3$ is a local maximum. If we denote the greatest value as 2, the least value as 0 and the value in the middle as 1, then we have a triplet (0,2,1). Any three consequitive numbers can be represented by such triplet. 
Out of six possible permutations\\
\hspace*{\fill}(0,1,2), (0,2,1), (1,0,2), (1,2,0), (2,0,1) and (2,1,0) \hspace*{\fill}\\
we are interested only in two combinations, which represent local maxima\\
\hspace*{\fill}(0,2,1) and (1,2,0).\hspace*{\fill}\\
Therefore, if we take any three consequitive numbers, then the probability that the number in the middle is a local maximum is 
	\[P_{max}=\frac{2}{6}=\frac{1}{3}\] 
	If we have $3\cdot N$ numbers in the sequence, then $N$ of them are local maxima. It also means that the average distance between maxima should be exactly $3$.
	
Now, let's introduce some additional notation and use an operator approach. When $x_{i-1} < x_i$ we'll put an operator $U$, i.e. the sequence goes "up". Alternatively, when $x_{i-1} > x_i$ we'll put an operator $D$, i.e. the sequence goes "down". For our sample sequence, the corresponding operator sequence is:
\begin{enumerate}
	\item 0.935536495
\item 0.191531578 D
\item 0.429049655 U
\item 0.308968021 D
\item 0.179540986 D
\item 0.401329789 U
\item 0.71581906 U
\item 0.604617962 D
\item 0.877254876 U
\item 0.973280207 U
\item 0.489033299 D
\item 0.912367351 U
\item 0.604552972 D
\item 0.039395302	D
\item 0.3780448	U
\item 0.55317569 U
\item 0.6308772	U
\item 0.373163479	D
\item 0.812434426 U
\item 0.560173882 D
\end{enumerate}

We can apply our new notation to triplets and see that\\ 
\hspace*{\fill}(0,1,2) becomes $\left\langle UU\right\rangle$ \hspace*{\fill}\\
\hspace*{\fill}(0,2,1) and (1,2,0) become $\left\langle UD\right\rangle$ \hspace*{\fill}\\
\hspace*{\fill}(1,0,2) and (2,0,1) become $\left\langle DU\right\rangle$ \hspace*{\fill}\\
\hspace*{\fill}(2,1,0) becomes $\left\langle DD\right\rangle$ \hspace*{\fill}\\
Having a new notation, we can use it to compute the probability of appearance of a 
local maximum in the middle of any triplet. We are interested in triplets represented by 
$\left\langle UD\right\rangle$, because it's the only expression, which represents a local maximum in the middle of a triplet. 
We shall use a standard cumulative distribution function (CDF) $F(x)$, defined as a probability of $x\leq x_{i}$: \[F(x_{i})=Pr(x\leq x_{i}) = \int^{x_i}_{-\infty}f(x)\cdot dx\] 
, where $f(x)$ is PDF (probability density function) of $x$. This also could be written as 
\[F(x_{i}) = \int^{F(x_i)}_{0}dF(x)\]
On the other hand, probability of $x>x_{i}$ is \[\int^{1}_{F(x_i)}dF(x)=1-F(x_{i})\] 
Now, we can write the following formula 
\[ 
\left\langle UD\right\rangle = \int^{1}_{0}dF(x_{i-1})\cdot \int^{1}_{F(x_{i-1})}dF(x_{i})\cdot  \int^{F(x_{i})}_{0}dF(x_{i+1})\]
The first integral declares that the first number in a triplet $\left\langle UD\right\rangle$ can have any value. The next integral says that the number in the middle of a triplet should be greater than the first number. Finally, the third integral is for the trailing number of a triplet, which should be less than the previous number. It's easy to compute the probability as follows
\[ 
\int^{1}_{0}dF(x_{i-1})\cdot \int^{1}_{F(x_{i-1})}dF(x_{i})\cdot  \int^{F(x_{i})}_{0}dF(x_{i+1}) =
\int^{1}_{0}dF(x_{i-1})\cdot \int^{1}_{F(x_{i-1})}dF(x_{i})\cdot F(x_{i})\]

\[ 
\int^{1}_{0}dF(x_{i-1})\cdot \int^{1}_{F(x_{i-1})}dF(x_{i})\cdot F(x_{i})=
\int^{1}_{0}dF(x_{i-1})\cdot ( \frac{1}{2}- \frac{F(x_{i-1})^2}{2})\]

\[ 
\int^{1}_{0}dF(x_{i-1})\cdot ( \frac{1}{2}- \frac{F(x_{i-1})^2}{2})=
\frac{1}{2}-\frac{1}{6}=\frac{1}{3}\]

We got the same number in both combinatorial and operator approaches. 

\section{PMF of distances between local maxima}

Now, it's time to advance our notation. Let's break up $\left\langle UD\right\rangle$ into pieces.
\begin{definition}
Operators $\langle$,$U$,$D$ and $\rangle$ are defined as

\[ 
\langle \psi(x) = \int^{1}_{0}d\psi(x)
\]

\[ 
U\cdot \psi(x) =  \int^{1}_{z(x)}d\psi(x)\cdot \psi(x)\]

\[ 
D\cdot \psi(x) =  \int^{z(x)}_{0}d\psi(x)\cdot \psi(x)\]

\[ 
\rangle = 1 
\]

\end{definition}

Armed with this notation let's look at any quintet of numbers from $i-1$ to $i+3$. 
If it happens so that the numbers come like (0,2,1,4,3), then we got  
two local maxima $x_{i}=2$ and $x_{i+2}=4$. The distance between these maxima is 
$(i+2) - i=2$. This quintet can be represented by an expression $\langle UDUD\rangle$. Such quintet in our sample sequence can be found at $i=10$: 
($x_9=0.877254876$, $x_{10}=0.973280207$, $x_{11}=0.489033299$,
$x_{12}=0.912367351$, $x_{13}=0.604552972$).

\begin{definition}
If $x_i$ is a local maximum, and the next nearest maximum is the number $x_j$, the the distance between maxima is $j-i$.
\end{definition}

Now, we can compute the probability of the distance between local maxima equal to 2
\begin{equation}
Pr(d=2)=f_m(2)=\frac{\langle UDUD\rangle }{P_{max}}=\frac{\langle UDUD\rangle }{\langle UD\rangle}
\end{equation}
, where $d$ is a distance between local maxima and $f_m(d)$ is a probability mass function (PMF) of the distribution of these distances.

Notice the denominator. It is necessary to divide the probability of the quintet by the probability of the maximum in its first three numbers (triplet). Consider the quintet ($x_{10}=0.973280207$,$x_{11}=0.489033299$,
$x_{12}=0.912367351$,$x_{13}=0.604552972$,$x_{14}=0.039395302$) from the sample sequence above. Its operator expression is $\langle DUDD\rangle$, which doesn't seem to represent two maxima on distance 2. Let's add the number $x_9=0.877254876$ and look at the resulting sextet. This sextet's operator expression is $\langle UDUDD\rangle$. It starts with $\langle UDUD\rangle$. Clearly, our original quintet is a part of a sextet with two maxima on distance 2. Therefore, we have to take into account those quintets, which are not represented by $\langle UDUD\rangle$, but these quintets could be parts of sought quintets. Probability $P_{max}$ in denominator includes those quintets, which would be left unaccounted otherwise. 

This formula can also be interpreted in terms of conditional probabilities as follows

\begin{equation}
Pr(A|B)=\frac{Pr(A\cap B) }{Pr(B)}
\end{equation}, where event $A\cap B$ is a quintet with one maximum in its head and one maximum in its tail, event B is the maximum in first three numbers of a quintet, and event $A|B$ is two maxima on a given distance from each other.

Using the same methodology it's easy to show that the probability of the distance 3 is
\begin{equation}
Pr(d=3)=f_m(3)=\frac{\langle UDUUD\rangle+\langle UDDUD\rangle }{P_{max}}=
\frac{\langle UDU^2D\rangle+\langle UD^2UD\rangle }{P_{max}}
\end{equation}
In order to see why there are two terms in the numerator, consider these two sextets (0,2,1,3,5,4) and (0,3,2,1,5,4). The following is the table with formulae for the next 3 distances

\begin{equation}
Pr(d=4)=f_m(4)=\frac{\langle UDU^3D\rangle+\langle UD^2U^2D\rangle+\langle UD^3D\rangle }
{P_{max}}
\label{pr4}
\end{equation}

\begin{equation}
Pr(d=5)=f_m(5)=\frac{\langle UDU^4D\rangle+\langle UD^2U^3D\rangle+\langle UD^3U^2D\rangle+\langle UD^4UD\rangle }
{P_{max}}
\end{equation}

\begin{equation}
f_m(6)=\frac{\langle UDU^5D\rangle+\langle UD^2U^4D\rangle+\langle UD^3U^3D\rangle+\langle UD^4U^2D\rangle+\langle UD^5UD\rangle }
{P_{max}}
\end{equation}

\section{Results}

In order to compute a probability of a given distance between maxima, we have to identify corresponding integrals, evaluate them and sum them up. 
For example, computing the probability of distance 4 involves evaluation, see equation \ref{pr4}. 
A simple analytical expression for the sums of integrals in numerators of the above probabilities was presented in 
\cite{oshanin} (see equation 3.8):
\begin{equation}
	p(l)=2^l\frac{(l-1)(l+2)}{(l+3)!}
	\label{pl}
\end{equation}
In \cite{oshanin} a set of similar problems are studied, e.g. permutation generated random walks, by using a different and more generic
approach. However, the equation \ref{pl} can be used to
to derive the probability of distances between local maxima:
\begin{equation}
	f_m(d)=\frac{p(l)}{P_{max}}=3\cdot 2^d\frac{(d-1)(d+2)}{(d+3)!}
	\label{pmf}
\end{equation}

We were not aware of this work, and in absence of a simple analytical expression for a sum of integrals in the PMF equations (such as \ref{pr4}), we wrote a Java program, which does all required work. First, it generates the necessary integrals using our operator notaion, e.g. $\langle UDU^3D\rangle+\langle UD^2U^2D\rangle+\langle UD^3D\rangle$. Next, it evaluates the corresponding integrals and sums symbolically. 

The resulting PMF table is shown in table \ref{pmf_table}. Variance of this distribution $\approx 1.167168$ and the standard deviation $\approx 1.08$. 

\begin{table}[h]
  \caption{Table of PMF of distances between local maxima}
  \label{pmf_table}
  \begin{center}
    \begin{tabular}{lll}
      Distance &Probability &Decimal Approximation\\
      $d$&$f_m(d)$\\
    \hline
2& 2/5 & 0.4\\
3& 1/3 & 0.3333333333333333\\
4& 6/35 & 0.17142857142857143\\
5& 1/15 & 0.06666666666666667\\
6& 4/189 & 0.021164021164021163\\
7& 1/175 & 0.005714285714285714 \\
8& 2/1485 & 0.0013468013468013469\\
9& 4/14175 & 2.821869488536155E-4\\
10& 4/75075 & 5.328005328005328E-5\\
11& 2/218295 & 9.161913923818686E-6\\
12& 4/2764125 & 1.4471125582236693E-6\\
13& 1/4729725 & 2.114287828573543E-7\\
14& 8/278326125 & 2.8743259368842937E-8\\
15& 1/273648375 & 3.654324641978963E-9\\
16& 2/4583103525 & 4.363855167334454E-10\\
17& 8/162820783125 & 4.913377669887681E-11\\
18& 4/764299911375 & 5.2335476433640715E-12\\
19& 2/3781060408125 & 5.289521414950853E-13\\
20& 4/78642438249375 & 5.0863122876684074E-14\\
21& 2/428772250281375 & 4.664480965565126E-15\\
22& 8/19566987612046875 & 4.088518968077954E-16\\
23& 2/58274046742786875 & 3.432059573325511E-17\\
24& 4/1447106344699640625 & 2.7641368684830376E-18\\
25& 8/37392513326621578125 & 2.1394657080470717E-19\\
26& 8/501914364595623354375 & 1.5938973985025013E-20\\
27& 4/3494761822449632109375 & 1.1445701318770342E-21\\
28& 8/100847608441898396203125 & 7.932761246003232E-23\\
29& 1/188217886723358757890625 & 5.312991328341671E-24\\
\cline{1-3}
\\
 Total&$\frac{2722885427931256697484374}{2722885427931256697484375 }$& 1 - 3.6725746509274224E-25

    \end{tabular}
    \par\medskip\footnotesize
  \end{center}
\end{table}

\begin{table}[h]
  \caption{Table of CDF of distances between local maxima}
  \label{cdf_table}
  \begin{center}
    \begin{tabular}{lll}
      Distance&Cumulative Probability &Decimal Approximation\\
      $d$&$F_m(d)$\\
    \hline
2&2/5&0.4\\
3&11/15&0.7333333333333333\\
4&19/21&0.9047619047619048\\
5&34/35&0.9714285714285714\\
6&134/135&0.9925925925925926\\
7&4717/4725&0.9983068783068783\\
8&5773/5775&0.9996536796536797\\
9&31183/31185&0.9999358666025333\\
10&184273/184275&0.9999891466558133\\
11&4729717/4729725&0.9999983085697371\\
12&16372121/16372125&0.9999997556822954\\
13&30405374/30405375&0.9999999671110782\\
14&241215974/241215975&0.9999999958543376\\
15&32564156609/32564156625&0.9999999995086623\\
16&36395233873/36395233875&0.9999999999450477\\
17&343732764373/343732764375&0.9999999999941815\\
18&3419236445623/3419236445625&0.999999999999415\\
19&142924083427117/142924083427125&0.999999999999944\\
20&782679504481871/782679504481875&0.9999999999999949\\
21&4482618980214373/4482618980214375&0.9999999999999996\\
22&53596531285171873/53596531285171875&$\approx 1.0$\\
23&$\frac{5341787618088796859}{5341787618088796875}$\\
24&$\frac{17307391882607701871}{17307391882607701875}$\\
25&$\frac{232984121496642140621}{232984121496642140625}$\\
26&$\frac{3253148659416077296871}{3253148659416077296875}$\\
27&$\frac{188217886723358757890609}{188217886723358757890625}$\\
28&$\frac{1408389014447201740078117}{1408389014447201740078125}$\\
29&$\frac{2722885427931256697484374}{2722885427931256697484375}$\\
\cline{1-3}
    \end{tabular}
    \par\medskip\footnotesize
  \end{center}
\end{table}

We tested validity of a computed PDF table on several random and pseudo-random number sequences. For pseudo-random number sequences we used Java's standard pseudo-random generator java.lang.Random and Daniel Cer's Java implementation \cite{Cer} of notorious RANDU generator \cite{Knuth}. 
We used Mads Haahr's True Random Number Service
web site \cite{Haahr} as a source of "true" random numbers. We modified the supplied Java client, which connects to the server and retrieves the true random number batches.  We generates random number sequences using these methods and compared them with the theoretical PMF using several tests such as Kolmogorov-Smirnov and $\chi^{2}$ goodness of fit tests, see chapters 1.3.5.15 and 1.3.5.16 in \cite{NIST}. Also, according to the central limit theorem in large samples the standard deviation of the average distance between maxima should approach $\frac {\sigma}{\sqrt{n}}$, where $\sigma$ is the standard deviation of the distances in the population and $n$ is the size of the sample \cite{MathWorld}. We used this feature to compare sample average distances to a theoretical average distance 3.

As expected, java.lang.Random's and "true" random sequences were consistent with our PMF on any sample sizes varying from 100 to 100,000,000. Surprisingly, RANDU-generated sequences were also compliant with this PMF. When deriving this PMF, we assumed that numbers in the sequences are independent of each other. RANDU generator has a well known deficiency: its numbers are not independent. However, as it was noted before, it fared well in our tests.

The table \ref{samples_tbl} shows sample statistics for "true" random and RANDU generated sequences compared to theoretical frequencies of distances between maxima. Both samples are distributed as predicted by theoretical PMF $f_m(d)$, they pass $\chi^{2}$ goodness of fit test with higher than 0.99 probabilities. Their average distances are also within the $3\cdot\sigma$ area of a theoretical value of 3. The p-value for the latter test is the probability of the deviation from the theoretical average distance greater than of the observed value.
\begin{table}[h]
	\caption{Sample Frequency Comparison}
	\label{samples_tbl}
	\centering
		\begin{tabular}{lrrr}
Distance&Theoretical&True Random&RANDU\\
&Frequency&Frequency&Frequency\\
\hline
2&40000&39803&40462\\
3&33333&33544&33003\\
4&17143&17119&17073\\
5&6667&6673&6545\\
6&2116&2139&2157\\
7&571&549&571\\
8&135&136&148\\
9&28&31&36\\
10&5&3&5\\
11&1&2&0\\
12&0&1&0\\
\hline
Average&3&3.00187&2.99447\\
Std Dev of Mean&0.0034\\
p-value&&0.584 &0.106\\
\hline
$\chi^{2}$, df&&1.006, 10&1.386, 8\\
p-value&&0.9998&0.9944\\
\cline{1-4}
		\end{tabular}
\end{table}

\begin{conclusion}
We constructed a simple method of computation of PMF of the distribution of distances between local maxima in random number series.
We confirmed that selected pseudo-random and true random number sequences are distributed according to this PMF.
\end{conclusion}
\begin{acknowledgement}
Author is very gratefull to Dr. F.M.Pen'kov and Dr. K.Gartvig for fruitfull discussions and interesting suggestions on this topic.
\end{acknowledgement}

\end{document}